\documentclass[10pt]{amsart}
\usepackage[T1]{fontenc}
\usepackage{graphicx}
\usepackage{amsthm,amsmath,amssymb,amsfonts,amscd}
\usepackage{graphicx}
\usepackage[all]{xy}
\usepackage{booktabs}
\usepackage{tikz}
\usepackage{subcaption}
\usepackage{xcolor,paralist,hyperref,fancyhdr,etoolbox}

\newcommand{\R}{\mathbb{R}}
\newcommand{\Ha}{\mathcal{H}}
\DeclareMathOperator{\diam}{diam}
\DeclareMathOperator{\dist}{dist}
\DeclareMathOperator{\cl}{cl}

\newcommand{\ldd}{\Theta_{\ast}^d}
\newtheorem*{theorem}{Theorem}
\newcommand{\vsp}{\vspace{10pt}}

\hypersetup{ colorlinks=true, linkcolor=black, filecolor=black, urlcolor=black }

\title[Besicovitch's example in higher dimensions]{Besicovitch's example in higher dimensions:\\a purely unrectifiable set with large lower density} 
\author{Jaume Capdevila Jové}
\date{}
\thanks{Partially supported by MICIU (Spain) under the grant PID2024-160507NB-I00.}
\address{Departament de Matemàtiques,
Universitat Autònoma de Barcelona, Bellaterra (Barcelona), Catalonia}
\email{jaume.capdevila.jove@uab.cat}

\begin{document}

\begin{abstract}
We generalize to arbitrary dimensions an example originally introduced by Besicovitch, obtaining for every $d \geq 1$ a purely $d$-unrectifiable set $E \subset \R^{d+1}$ such that $\Theta_{\ast}^{d}(E, x) = \liminf_{r \to 0} \Ha^{d}(E \cap B_{r}(x))/(2r)^{d} = 1/2$ for $\Ha^{d}$-almost every point $x \in E$. This establishes the lower bound $1/2$ for the minimal value $\sigma$ such that, if $\Theta_{\ast}^{d}(E, x) > \sigma$ for $\Ha^{d}$-almost all $x \in E$, then $E$ is $d$-rectifiable. This threshold was conjectured to be exactly $1/2$ by Besicovitch.
\end{abstract}

\maketitle


\section{Introduction}

One of the main concepts of geometric measure theory is that of $d$-rectifiable subsets of $\mathbb{R}^{n}$, given positive integers $0 < d < n$. They are sets which, up to a set of zero $\mathcal{H}^{d}$-measure, are contained in a countable union of images of Lipschitz maps with domain in $\mathbb{R}^{d}$ (where $\mathcal{H}^{d}$ denotes the $d$-dimensional Hausdorff measure). For example, for $d = 1$, the $1$-rectifiable sets are those which are contained in a countable union of rectifiable curves, again up to a set of zero $\mathcal{H}^{1}$-measure. On the other side of the coin, we have the purely $d$-unrectifiable sets, which are those that contain no $d$-rectifiable subset of positive $\mathcal{H}^{d}$-measure. One of the goals of geometric measure theory is to characterize rectifiability in terms of other geometric or analytical properties.

To that end, one of the basic tools is that of the densities for the Hausdorff measure. Consider a Borel set $E \subset \mathbb{R}^{n}$ such that $0 < \mathcal{H}^{d}(E) < \infty$ for some positive integers $0 < d < n$, which we call a $d$-set. One defines the upper and lower $d$-densities of $E$ at a point $x \in \mathbb{R}^{n}$, denoted as $\Theta^{\ast d}(E,x)$ and $\Theta_{\ast}^{d}(E, x)$ respectively, as the $\limsup$ and $\liminf$ as $r \to 0$ of
\begin{equation*}
    \frac{\mathcal{H}^{d}(E \cap B_{r}(x))}{ (2r)^{d}}.
\end{equation*}
When both quantities coincide, the limit is called the $d$-density of $E$ at $x$. A fundamental theorem states that a $d$-set $E \subset \mathbb{R}^{n}$ is $d$-rectifiable if and only if the $d$-density of $E$ exists and is equal to $1$ at $\mathcal{H}^{d}$-almost all points of $E$. This is known as the characterization of rectifiability in terms of densities. This line of study was initiated in the pioneering work of Besicovitch \cite{besicovitch1938} in 1938, where he established the result for $1$-sets in the plane, i.e., the case $d = 1$ and $n = 2$. It was extended to arbitrary dimension over different stages, with the work of Moore \cite{moore} (1950), Marstrand \cite{marstrand} (1961) and Mattila \cite{mattila-2} (1975). It was later proven by Preiss \cite{preiss} in 1987 that in fact, rectifiability is equivalent to the mere existence of the density at almost all points.

In his original article, Besicovitch also proved that if $\Theta_{\ast}^{1}(E, x) > 3/4$ for $\mathcal{H}^{1}$-almost all points of a $1$-set $E$, then $E$ is automatically $1$-rectifiable. Following this idea, we may naturally define the following coefficient:
\begin{equation*}
\begin{split}
\sigma_{d,n} := \min\{&\sigma > 0\ :\ \text{for any $d$-set $E\subset \mathbb{R}^{n}$,}\\
&\Theta_{\ast}^{d}(E,x) > \sigma\ \mathcal{H}^{d}\text{-a.e. }x\in{E} \implies E\text{ is $d$-rectifiable.}\}
\end{split}
\end{equation*}

Equivalently, $\sigma_{d,n}$ is defined as the greatest value that the lower $d$-density of a purely $d$-unrectifiable set can attain at $\Ha^{d}$-almost all of its points. The previously stated result of Besicovitch translates to the bound $\sigma_{1,2} \leq 3/4$. Moreover, in the same article in 1938 he provided an example of a purely $1$-unrectifiable set $P$ contained in the plane which satisfies $\Theta^{1}_{\ast}(P, x) = 1/2$ at $\mathcal{H}^{1}$-almost all $x \in P$; a formal proof of this fact appeared later in a paper by Dickinson \cite{dickinson} in 1939. This way, they proved the lower bound $\sigma_{1,2} \geq 1/2$. With this in mind, Besicovitch conjectured that the exact value of $\sigma_{1,2}$ is $1/2$, which is now known as \emph{Besicovitch's $1/2$-conjecture}, and is still an open problem. This conjecture extends to higher dimensions, and the value of $\sigma_{d, n}$ is expected to be $1/2$ for all $0 < d < n$.

Further improvements to this bound have been obtained since then. In 1992, Preiss and Tišer \cite{preiss-tiser} refined the estimate to $\sigma_{1,n} \leq (2+\sqrt{46})/12 < 59/80$, which holds for all $n \geq 2$ (for all metric spaces, in fact). Recently, in 2024, Camillo De Lellis et al. \cite{camilo-lellis} established that $\sigma_{1,n} \leq 7/10$, which is currently the best known upper bound.

In higher dimensions (for $d > 1$), no good upper bounds are known for $\sigma_{d,n}$. On the other hand, the same lower bound remains valid; in this work, we generalize Besicovitch’s example to arbitrary dimensions, thereby proving the following theorem.
\begin{theorem}
For all $0 < d < n$,
\begin{equation*}
    \sigma_{d,n} \geq 1/2.
\end{equation*}
\end{theorem}
As far as we know, this generalization has not appeared in the literature before. We will construct the set in codimension one, hence work in $n = d+1$. The trivial inclusion of this set into arbitrary $n > d+1$ yields the general result.

\section{Besicovitch's example}

We begin by constructing a set $\Pi \subset \R^{2}$, which is the example in the plane originally introduced by Besicovitch. Let $\alpha > 0$, $\beta > 0$, and consider the function
\begin{equation*}
    f(\alpha, \beta, x) = \begin{cases}
        -\beta\quad &\text{if}\quad 2k\alpha \leq x < (2k+1)\alpha\\
        \beta\quad &\text{if}\quad (2k+1)\alpha \leq x < 2(k+1)\alpha,
    \end{cases}
\end{equation*}
for any integer $k$. Given
\begin{equation*}
    \alpha_{n} = 2^{-n^{2}},\quad \beta_{n} = \frac{1}{n}2^{-n^{2}},
\end{equation*}
we consider $f_{n}(x) = f(\alpha_{n}, \beta_{n}, x)$, and the sum $g(x) = \sum_{n = 1}^{\infty}f_{n}(x)$. We define
\begin{equation*}
    \Pi = \left\{(x, g(x))\ :\ x \in [0,1]\right\}
\end{equation*}
as the graph of $g$ above the interval $[0, 1]$. For each positive integer $n$, we consider the sets of points
\begin{equation*}
    \mathcal{C}_{n} = \{ i\cdot 2^{-n^{2}}:0 \leq i \leq 2^{n^{2}}\},
\end{equation*}
and the intervals
\begin{equation*}
    I_{n,i} = [i\cdot \alpha_{n}, (i+1)\alpha_{n}), \quad \text{for $0 \leq i < 2^{n^{2}}$}.
\end{equation*}
Notice that if $x, y \in I_{n,i}$, then $f_{j}(x) = f_{j}(y)$ for all $1 \leq j \leq n$, and so
\begin{equation}
\label{EQ:example upper bound}
\begin{split}
    |g(x) - g(y)| &\leq \sum_{j = 1}^{\infty}\left|f_{j}(x) - f_{j}(y)\right| =  \sum_{j = n+1}^{\infty} \left|f_{j}(x) - f_{j}(y)\right|\\
    & \leq 2\sum_{j = n+1}^{\infty}\frac{1}{j}2^{-j^{2}} \leq \frac{2}{n} \sum_{j = (n+1)^{2}}^{\infty}2^{-j} \leq \frac{4}{n}2^{-(n+1)^{2}}  = \alpha_{n+1}\frac{4}{n}.
\end{split}
\end{equation}
On the other hand, if $x\in I_{n,i}$ and $y \in I_{n,i-1}\cup I_{n,i+1}$ then there is some $k \leq n$ such that $f_{k}(x) \neq f_{k}(y)$ whilst $f_{j}(x) = f_{j}(y)$ for all $1 \leq j < k$, so
\begin{equation}
\label{EQ:example lower bound}
\begin{split}
    |g(x) - g(y)| & = \left|\sum_{j = 1}^{\infty}f_{j}(x) - \sum_{j = 1}^{\infty} f_{j}(y)\right| \\
    & = \left|f_{k}(x) - f_{k}(y) + \sum_{j = k+1}^{\infty} (f_{j}(x) - f_{j}(y))\right|\\
    & \geq \left||f_{k}(x) - f_{k}(y)| - \left|\sum_{j = k+1}^{\infty}(f_{j}(x) - f_{j}(y))\right|\right|  \geq \frac{2\alpha_{k}}{k} - \frac{4\alpha_{k+1}}{k} \geq \frac{\alpha_{n}}{n}.\\
\end{split}
\end{equation}

From now on, we fix the dimension $d \geq 1$. We will show that the Cartesian product
\begin{equation*}
    \Pi_{d} = \Pi \times [0,1]^{d-1} \subset \R^{d+1}
\end{equation*}
is a purely $d$-unrectifiable $d$-set that has lower $d$-density equal to $1/2$ $\Ha^{d}$-almost everywhere. We define the ($d$-dimensional) rectangles $R_{n,i} = I_{n,i}\times [0,1]^{d-1} \subset \R^{d}$. Moreover, we will write $d_{n}(x) = \dist(x_{1},\mathcal{C}_{n})$, where now we denote $x = (x_{1},\dots, x_{d}) \in \R^{d}$. Notice that up to a reordering of the coordinates,
\begin{equation*}
    \Pi_{d} = \left\{(x_{1},\dots, x_{d}, g(x_{1}))\ :\ (x_{1},\dots, x_{d}) \in [0,1]^{d}\right\},
\end{equation*}
which can also be understood as the graph of the function $G:[0,1]^{d} \to \R$ defined by $G(x) = g(x_{1})$. As before, for any $x \in [0,1]^{d}$ we denote $\Pi_{d}(x) = (x, G(x)) \in \Pi_{d}$, so for any subset $E \subset [0,1]^{d}$, the set $\Pi_{d}(E) \subset \Pi_{d}$ is the part of the graph above $E$. In particular, following this notation we have $\Pi_{d} = \Pi_{d}([0,1]^{d})$.  Finally, observe that for the function $G$, Eq. \ref{EQ:example upper bound} holds similarly if $x,y\in R_{n,i}$, and Eq. \ref{EQ:example lower bound} holds similarly if $x\in R_{n,i}$ and $y \in R_{n,i-1}\cup R_{n,i+1}$.

\vsp

We divide the proof into four parts, \textbf{(A)}, \textbf{(B)}, \textbf{(C)} and \textbf{(D)}.

\vsp

Before starting, we clarify that we use the convention
\begin{equation*}
    \Ha^{d}_{\delta}(E) = \inf \left\{\sum_{i = 1}^{\infty}\diam(U_{i})^{d}\ :\ E \subset \bigcup_{i = 1}^{\infty}U_{i}, \ \diam(U_{i})\leq \delta\right\},
\end{equation*}
and $\Ha^{d}(E) = \lim_{\delta \to 0}\Ha^{d}_{\delta}(E)$ as usual. This way, for any ball $B$ we have that $\Ha^{d}(B) = \diam(B)^{d}$, and for any cube $Q$ we have that $\Ha^{d}(Q) = c_{d} \diam(Q)^{d}$, where $c_{d}$ is a constant depending only on $d$.

\vsp
\noindent \textbf{(A)} For any Borel set $Z \subset[0,1]^{d}$ such that $\Ha^{d}(Z) = 0$, we have $\Ha^{d}(\Pi_{d}(Z)) = 0$.
\vsp

For every $n \geq 1$, we partition $[0,1]^{d}$ (except for a subset of the boundary, which has zero $\Ha^{d}$-measure) into a union of ($d$-dimensional) cubes of side length $\alpha_{n}$. Precisely, $Q_{n,i} = I_{n,i_{1}}\times\dots\times I_{n,i_{d}}$ where $i = (i_{1},\dots, i_{d})$ and $0 \leq i_{j} < 2^{n^{2}}$ for any $1 \leq j \leq d$. We denote by $\mathcal{Q}_{n}$ the set of such cubes $Q_{n,i}$ of generation $n$. Due to Eq. \ref{EQ:example upper bound} we can cover $\Pi_{d}(Q_{n,i})$ with a ($(d+1)$-dimensional) rectangle $\tilde{Q}_{n,i}$ that is the Cartesian product of $Q_{n,i}$ with an interval of length $\alpha_{n}/n$. Then $\diam(Q_{n,i}) = \alpha_{n}\sqrt{d}$ and
\begin{equation*}
    \diam(\tilde{Q}_{n,i}) = \sqrt{\diam(Q_{n,i})^{2} + \left(\frac{\alpha_{n}}{n}\right)^{2}} = \diam(Q_{n,i})\sqrt{1+\frac{1}{dn^{2}}}.
\end{equation*}
Now, let $E \subset [0,1]^{d}$ be a compact set. Denote by $\mathcal{Q}_{n}(E)$ the set of cubes $Q \in \mathcal{Q}_{n}$ such that $Q\cap E \neq \varnothing$. Then $\Pi_{d}(E) \subset \cup_{Q \in \mathcal{Q}_{n}(E)} \tilde{Q}$, and setting $\delta_{n} = \alpha_{n}\sqrt{d + 1/n^{2}}$, we have
\begin{equation}
\label{EQ:example 1}
    \Ha^{d}_{\delta_{n}}(\Pi_{d}(E)) \leq \sum_{Q \in \mathcal{Q}_{n}(E)}\diam(\tilde{Q})^{d} = \left(1+\frac{1}{dn^{2}}\right)^{\frac{d}{2}}\sum_{Q \in \mathcal{Q}_{n}(E)}\diam(Q)^{d}.
\end{equation}
Since $E$ is compact, we can write
\begin{equation*}
    E = \bigcap_{n = 1}^{\infty}\bigcup_{Q \in \mathcal{Q}_{n}(E)}\cl(Q),
\end{equation*}
where $\cl(Q)$ denotes the closure of $Q$, and in particular
\begin{equation*}
    \Ha^{d}(E) = \lim_{n \to \infty} \sum_{Q \in \mathcal{Q}_{n}(E)}\Ha^{d}(Q) = c_{d}\lim_{n \to \infty} \sum_{Q \in \mathcal{Q}_{n}(E)}\diam(Q)^{d}.
\end{equation*}
Taking limits as $n \to \infty$ of Eq. \ref{EQ:example 1}, we get
\begin{equation}
\label{EQ:example aa}
    \Ha^{d}(\Pi_{d}(E)) \leq c_{d}^{-1}\Ha^{d}(E).
\end{equation}

Now, let $Z\subset [0,1]^{d}$ be a general Borel set such that $\Ha^{d}(Z) = 0$. We argue by contradiction and suppose that $\Ha^{d}(\Pi_{d}(Z)) > 0$. By the regularity of Hausdorff measure (see Theorem 1.6 of Falconer \cite{falconer}), we can find a compact subset $X \subset \Pi_{d}(Z)$  such that $\Ha^{d}(X) > 0$. But $X = \Pi_{d}(E)$ for some $E \subset Z$, and $E$ is also compact because it is the image of the compact set $\Pi_{d}(E)$ under the projection into $[0,1]^{d}$. Using Eq. \ref{EQ:example aa}, $\Ha^{d}(X) = \Ha^{d}(\Pi_{d}(E)) \leq c_{d}^{-1}\Ha^{d}(E) \leq c_{d}^{-1} \Ha^{d}(Z) = 0$, which is a contradiction. Therefore, we must have $\Ha^{d}(\Pi_{d}(Z)) = 0$.

\vsp
\noindent \textbf{(B)} For any Borel set $E \subset [0,1]^{d}$, we have $\Ha^{d}(\Pi_{d}(E)) = \Ha^{d}(E)$.
\vsp

By basic properties of the Hausdorff measure under Lipschitz mappings (see Lemma 1.8 of Falconer \cite{falconer}), we obtain 
\begin{equation}
\label{EQ:example d projection}
\Ha^{d}(E) \leq \Ha^{d}(\Pi_{d}(E))
\end{equation}
by projecting into $[0,1]^{d}$. For the reverse inequality, it will suffice to show that $\Ha^{d}(\Pi_{d}) = \Ha^{d}([0,1]^{d})$. Indeed, if this is true, denoting $E' = [0,1]^{d}\setminus E$, we have that if $\Ha^{d}(\Pi_{d}(E)) > \Ha^{d}(E)$ holds, then
\begin{equation*}
\begin{split}
    \Ha^{d}(\Pi_{d}(E')) &= \Ha^{d}(\Pi_{d}) - \Ha^{d}(\Pi_{d}(E)) \\
    &< \Ha^{d}([0,1]^{d}) - \Ha^{d}(E) = \Ha^{d}(E'),
\end{split}
\end{equation*}
which contradicts Eq. \ref{EQ:example d projection}.

So, let us prove that $\Ha^{d}(\Pi_{d}) = \Ha^{d}([0,1]^{d})$. For that, we first need some constructions. For each cube $Q \in \mathcal{Q}_{n}$, we can consider the inscribed closed ball $B$, which has diameter equal to the side length of $Q$, $\diam(B) = \alpha_{n}$. For each $n \geq 1$, we define $\mathcal{V}_{n}$ as the set of all inscribed balls in cubes $Q \in \mathcal{Q}_{m}$ for all $m \geq n$. Each $B \in \mathcal{V}_{n}$ is contained in a cube $Q \in \mathcal{Q}_{m(B)}$ for some $m(B) \geq n$, such that $\diam(B) = \alpha_{m(B)}$. Hence $\Pi_{d}(B)$ is covered by the product of $B$ with an interval of length $\alpha_{m(B)}/m(B)$, that we denote $C(B)$, and which satisfies
\begin{equation}
\label{EQ:example cylinder}
\begin{split}
    \diam(C(B)) & = \sqrt{\diam(B)^{2}+\left(\frac{\alpha_{m(B)}}{m(B)}\right)^{2}} \\
    & = \diam(B)\sqrt{1 + \frac{1}{m(B)^{2}}} \leq \diam(B)\sqrt{1+\frac{1}{n^{2}}}.
\end{split}
\end{equation}
We will remove a set of zero $\Ha^{d}$-measure from $[0,1]^{d}$, so that every point in the remaining set is contained in arbitrarily small balls of $\mathcal{V}_{n}$. Given a cube $Q \subset \R^{d}$ of side length $2\ell$, of center $(x_{1},\dots, x_{d})$ and oriented parallel to the coordinate axis, we consider the inner region
\begin{equation*}
    \text{In}(Q) = \left\{(y_{1},\dots, y_{d}) \in \R^{d}\ :\ \sum_{i = 1}^{d}|y_{i}-x_{i}| \leq \ell\right\}.
\end{equation*}
This region is contained in the ball centered at $(x_{1}, \dots, x_{d})$ with radius $\ell$, and it has volume $2^{d}\ell^{d}/d!$. The ratio of $\Ha^{d}$-measure of the exterior region of $Q$, which we denote $\text{Ex}(Q) = Q\setminus \text{In}(Q)$, is then $\gamma_{d} = 1 - 1/d! < 1$. Now take a cube $Q \in \mathcal{Q}_{n}$, keep the exterior region $\text{Ex}(Q)$, and consider $\{Q'\cap \text{Ex}(Q)\ :\ Q' \in \mathcal{Q}_{n+1}\}$. For these new cubes (and ``half cubes'', as some new cubes are split apart when taking the intersection), we repeat the process indefinitely. For each $N \geq 1$, we define
\begin{equation*}
    \mathcal{X}_{N} = \bigcup_{Q \in \mathcal{Q}_{N}} \bigcap_{m \geq 1} \text{Ex}^{m}(Q),
\end{equation*}
where $\text{Ex}^{m}(Q)$ is the remaining set after applying the procedure $m$ times. We have $\Ha^{d}(\text{Ex}^{m}(Q)) = \gamma_{d}^{m}\Ha^{d}(Q) \to 0$ as $m \to \infty$, so $\Ha^{d}(\mathcal{X}_{N}) = 0$. Hence, if we define $\mathcal{X} = \cup_{N \geq 1} \mathcal{X}_{N}$ we also have $\Ha^{d}(\mathcal{X}) = 0$.

By construction, for any $n \geq 1$, every $x \in F := [0,1]^{d}\setminus \mathcal{X}$ is contained in arbitrarily small balls of $\mathcal{V}_{n}$. By Vitali's covering theorem (see Theorem 1.10 of Falconer \cite{falconer}), there is a countable collection of pairwise disjoint balls $\mathcal{B}_{n} \subset \mathcal{V}_{n}$ such that
\begin{equation*}
    \Ha^{d}\left(F \setminus \bigcup_{B \in \mathcal{B}_{n}}B\right) = 0.
\end{equation*}
In other words, there is a Borel set $Z_{n}$ with $\Ha^{d}(Z_{n}) = 0$ such that $F \setminus Z_{n} \subset \cup_{B \in \mathcal{B}_{n}}B$. Now consider $Z = \cup_{n \geq 1} Z_{n}$, which is also a Borel set with $\Ha^{d}(Z) = 0$, and take $\tilde{F} = [0,1]^{d} \setminus (\mathcal{X}\cup Z)$. By \textbf{(A)}, $\Ha^{d}(\Pi_{d}(\mathcal{X}\cup Z)) = 0$, so $\Ha^{d}(\Pi_{d}) = \Ha^{d}(\Pi_{d}(\tilde{F}))$. For each $n \geq 1$, consider the covering of $\Pi_{d}(\tilde{F})$ given by the cylinders $C(B)$ for $B \in \mathcal{B}_{n}$, and using Eq. \ref{EQ:example cylinder} we have
\begin{equation*}
    \begin{split}
        \Ha^{d}_{\delta_{n}}(\Pi_{d}(\tilde{F})) &\leq \sum_{B \in \mathcal{B}_{n}}\diam(C(B))^{d} \\
        & \leq \left(1 + \frac{1}{n^{2}}\right)^{\frac{d}{2}}\sum_{B \in \mathcal{B}_{n}}\Ha^{d}(B)\\
        & = \left(1 + \frac{1}{n^{2}}\right)^{\frac{d}{2}}\Ha^{d}([0,1]^{d}).
    \end{split}
\end{equation*}
Taking the limit $n \to \infty$, we find
\begin{equation*}
    \Ha^{d}(\Pi_{d}) = \Ha^{d}(\Pi_{d}(\tilde{F})) \leq \Ha^{d}([0,1]^{d}),
\end{equation*}
so finally
\begin{equation*}
    \Ha^{d}(\Pi_{d}) = \Ha^{d}([0,1]^{d}).
\end{equation*}

\vsp
\noindent \textbf{(C)} $\ldd(\Pi_{d}, \Pi_{d}(x)) \geq 1/2$ for all $x \in (0,1)^{d}$.
\vsp

\begin{figure}
\centering
\begin{tikzpicture}[scale=0.9, every node/.style={transform shape}]
    \begin{scope}
        \draw[thick] (0,0) circle(2);
        \begin{scope}
            \clip (0,0) circle(2);
            \fill[gray!30] (-2.5,0.8) rectangle (2.5,2.5);
        \end{scope}
        \draw[thick, black] (-2.5,0.8) -- (2.5,0.8);
        \filldraw (0,0) circle(1pt);
        \draw[thick] (2.6,0) --++ (0.2,0) 
                      --++ (0,0.8) 
                      --++ (-0.2,0);
        \node[right] at (2.8,0.4) {\(h\)};
        \draw[thick] (-2,0) --++ (0,-0.2) 
                            --++ (2,0) 
                            --++ (0,0.2);
        \node[below] at (-1, -0.3) {\(r\)};
        \node at (0, 1.2) {\(\text{SC}(r,h)\)};
        \draw[dashed] (-2.5,0) -- (2.5,0);
    \end{scope}

    \begin{scope}[shift={(7,0)}]
        \draw[thick] (0,0) circle(2);
        \draw[thick, black] (-2.5,0.75) -- (2.5,0.75);
        \draw[thick, black] (-2.5,-1.4) -- (2.5,-1.4);
        \begin{scope}
            \clip (0,0) circle(2);
            \fill[gray!30] (-2.5,-1.4) rectangle (2.5,0.75);
        \end{scope}
        \filldraw (0,0) circle(1pt);
        \draw[thick] (-2,0) --++ (0,0.2) 
                            --++ (2,0) 
                            --++ (0,-0.2);
        \node[above] at (-1, 0.3) {\(r\)};
        \draw[thick] (-2.6,0) --++ (-0.2,0) 
                              --++ (0,0.75) 
                              --++ (0.2,0);
        \node[left] at (-2.8,0.375) {\(a\)};
        \draw[thick] (2.6,-1.4) --++ (0.2,0) 
                                --++ (0,2.15) 
                                --++ (-0.2,0);
        \node[right] at (2.8,-0.325) {\(c\)};
        \draw[dashed] (-2.5,0) -- (2.5,0);
        \node at (0,-0.6) {\( A(r, a, c) \)};
    \end{scope}

\end{tikzpicture}
\caption{Picture for $d = 2$. On the left, a spherical cap $SC(r,h)$. On the right, the spherical region $A(r, a, c)$.}
\label{FIG:example 1}
\end{figure}
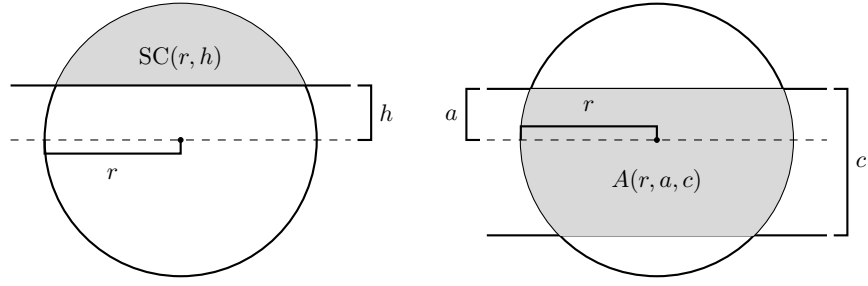

We need to introduce some notation for geometric objects concerning hyperspheres, which in the case $d = 1$ would simply be intervals. We denote by $SC(r, h)$ a $d$-dimensional spherical cap of radius $r$ at height $0 \leq h \leq r$ from the center (see Figure \ref{FIG:example 1}). Using Fubini's theorem, one can easily compute
\begin{equation}
\label{EQ:area SC}
        \Ha^{d}(SC(r, h)) = b_{d} r^{d} \int_{0}^{\arccos(h/r)}\sin^{d}(\theta)d\theta,
\end{equation}
where $b_{d}$ is a constant depending only on $d$. We also denote by $A(r, a, c)$ the region consisting of a $d$-dimensional ball minus two parallel $d$-dimensional spherical caps in opposite hemispheres, one at height $0\leq a \leq r$, and the other at height $0 \leq c - a\leq r$ (see Figure \ref{FIG:example 1}). These geometric objects are considered in $\R^{d}$. To ease the notation, we will omit the explicit dependence on the dimension for both these objects and the $d$-dimensional balls of $\R^{d}$. On the contrary, when considering $(d+1)$-dimensional balls in $\R^{d+1}$, we will denote it with a superscript to distinguish them from the rest.

Now, consider the quantities $\lambda_{n} = \sqrt{1+(4/n)^{2}}$, which are decreasing in $n$. Observe that for $n \geq 2$, $(\alpha_{n+1}\lambda_{n})/(\alpha_{n} \lambda_{n-1}) < 2^{-2n + 1} < 1$, so we always have $\alpha_{n+1}\lambda_{n} < \alpha_{n}\lambda_{n-1}$. Moreover, any radius in the interval $(0, 1/4)$ is contained in some interval $[\alpha_{n+1}\lambda_{n}, \alpha_{n}\lambda_{n-1}]$. So, let $x \in [0,1]^{d}$ and take any radius $0 < r < 1/4$. Then, let $n \geq 2$ such that
\begin{equation}
\label{EQ:example alpha intervals}
\alpha_{n+1}\lambda_{n} \leq r \leq \alpha_{n}\lambda_{n-1},
\end{equation}
and $0 \leq i < 2^{n^{2}}$ be such that $x \in R_{n,i}$. Then, for any $y \in R_{n,i} \cap B_{r/\lambda_{n}}(x)$, using $|x - y | \leq r/\lambda_{n}$ and Eq. \ref{EQ:example upper bound}, we find
\begin{equation*}
\begin{split}
    \dist\left(\Pi_{d}(x), \Pi_{d}(y)\right) &= \sqrt{\left|x-y\right|^2+\left|G(x)-G(y)\right|^{2}} \\
    & \leq \sqrt{\left(\frac{r}{\lambda_{n}}\right)^{2}+\left(\alpha_{n+1}\lambda_{n}\right)^{2}\left(\frac{4}{n\lambda_{n}}\right)^{2}} \\
    & \leq \frac{r}{\lambda_{n}}\sqrt{1 + \left(\frac{4}{n}\right)^{2}} = r,
\end{split}
\end{equation*}
where in the last inequality we have also used the lower bound of Eq. \ref{EQ:example alpha intervals}. Thus,
\begin{equation}
\label{EQ:generalization inclusion 1}
    \Pi_{d}\left(R_{n,i}\cap B_{r/\lambda_{n}}(x)\right) \subset \Pi_{d} \cap B_{r}^{d+1}(\Pi_{d}(x)).
\end{equation}

\noindent \emph{Case 1:} $r/\lambda_{n} \leq \alpha_{n}-d_{n}(x)$.

\begin{figure}[htbp]
\centering
\begin{subfigure}[b]{0.49\textwidth}
    \centering
    \begin{tikzpicture}[scale=0.7, every node/.style={transform shape}]
        \begin{scope}
            \draw[thick] (0,0) circle(2);
            \draw[thick, black] (-2.8,0.75) -- (2.8,0.75);
            \draw[thick, black] (-2.8,-2.2) -- (2.8,-2.2);
            \begin{scope}
                \clip (0,0) circle(2);
                \fill[gray!30] (-3,0.75) rectangle (3,-2.2);
            \end{scope}
            \filldraw (0,0) circle(1pt);
            \draw[thick] (-2,0) --++ (0,-0.2) 
                                 --++ (2,0) 
                                 --++ (0,0.2);
            \node[below] at (-1, -0.3) {$r/\lambda_{n}$};
            
            \draw[thick] (3,0) --++ (0.15,0) 
                               --++ (0,0.75) 
                               --++ (-0.15,0);
            \node[right] at (3.2,0.375) {$d_{n}(x)$};

            \draw[thick] (-3,-2.2) --++ (-0.15,0) 
                                   --++ (0,2.95) 
                                   --++ (0.15,0);
            \node[left] at (-3.2,-0.725) {$\alpha_{n}$};

            \draw[thick] (3,0) --++ (0.15,0) 
                               --++ (0,-2.2) 
                               --++ (-0.15,0);
            \node[right] at (3.2,-1.1) {$\alpha_{n} - d_{n}(x)$};

            \draw[dashed] (-2,0) -- (2,0);
            
            \node at (-2.2, 1.4) {$R_{n,i+1}$};
            \node at (-2.4, -0.75) {$R_{n,i}$};
            \node at (-2.2, -2.6) {$R_{n,i-1}$};
            \node at (0, 0.3) {$x$};
        \end{scope}
    \end{tikzpicture}
    \caption{Case $r/\lambda_{n} \leq \alpha_{n}-d_{n}(x)$.}
    \label{FIG:example 2}
\end{subfigure}
\hfill
\begin{subfigure}[b]{0.49\textwidth}
    \centering
    \begin{tikzpicture}[scale=0.7, every node/.style={transform shape}]
        \begin{scope}
            \draw[thick] (0,0) circle(2.5);
            \draw[thick, black] (-3.3,0.75) -- (3.3,0.75);
            \draw[thick, black] (-3.3,-2.2) -- (3.3,-2.2);
            \begin{scope}
                \clip (0,0) circle(2.5);
                \fill[gray!30] (-3.5,0.75) rectangle (3.5,-2.2);
            \end{scope}
            \filldraw (0,0) circle(1pt);
            \draw[thick] (-2.5,0) --++ (0,-0.2) 
                                 --++ (2.5,0) 
                                 --++ (0,0.2);
            \node[below] at (-1.25, -0.3) {$r/\lambda_{n}$};
            
            \draw[thick] (3.5,0) --++ (0.15,0) 
                                 --++ (0,0.75) 
                                 --++ (-0.15,0);
            \node[right] at (3.7,0.375) {$d_{n}(x)$};

            \draw[thick] (-3.5,-2.2) --++ (-0.15,0) 
                                   --++ (0,2.95) 
                                   --++ (0.15,0);
            \node[left] at (-3.7,-0.725) {$\alpha_{n}$};

            \draw[thick] (3.5,0) --++ (0.15,0) 
                                 --++ (0,-2.2) 
                                 --++ (-0.15,0);
            \node[right] at (3.7,-1.1) {$\alpha_{n} - d_{n}(x)$};

            \draw[dashed] (-2.5,0) -- (2.5,0);
            
            \node at (-2.7, 1.4) {$R_{n,i+1}$};
            \node at (-2.9, -0.75) {$R_{n,i}$};
            \node at (-2.7, -2.6) {$R_{n,i-1}$};
            \node at (0, 0.3) {$x$};
        \end{scope}
    \end{tikzpicture}
    \caption{Case $r/\lambda_{n} > \alpha_{n}-d_{n}(x)$.}
    \label{FIG:example 3}
\end{subfigure}

\caption{Picture for $d = 2$. The region $R_{n,i}\cap B_{r/\lambda_{n}}(x)$ is coloured in grey for both cases.}
\label{FIG:examples_combined}
\end{figure}
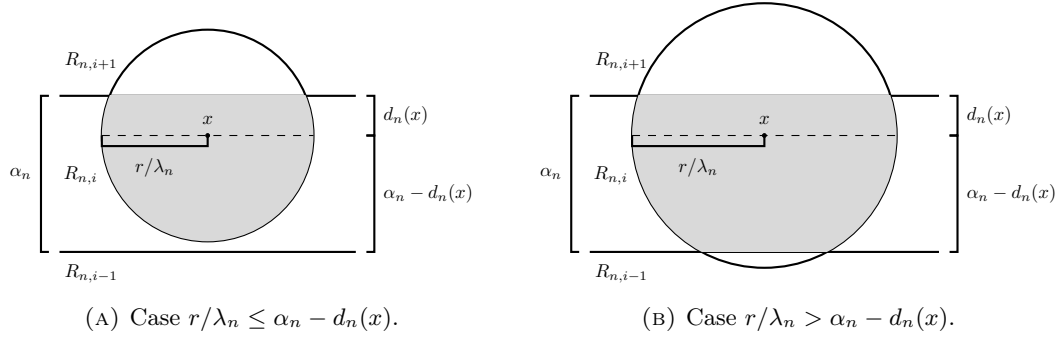

In this case, $R_{n,i}\cap B_{r/\lambda_{n}}(x)$ consists of $B_{r/\lambda_{n}}(x)$ without a spherical cap $SC(r/\lambda_{n}, d_{n}(x))$ (see Figure \ref{FIG:example 2}). In particular, it contains half of a $d$-dimensional ball of radius $r/\lambda_{n}$. Hence, using Eq. \ref{EQ:generalization inclusion 1} and part \textbf{(B)},


\begin{equation*}
\begin{split}
    \frac{\Ha^{d}\left(\Pi_{d} \cap B_{r}^{d+1}(\Pi_{d}(x))\right)}{(2r)^{d}} &\geq \frac{\Ha^{d}\left(\Pi_{d}\left(R_{n,i}\cap B_{r/\lambda_{n}}(x)\right)\right)}{(2r)^{d}}\\
    &= \frac{\Ha^{d}\left(R_{n,i} \cap B_{r/\lambda_{n}}(x)\right)} {(2r)^{d}}
    \\
    & \geq \frac{\frac{1}{2}(2r/\lambda_{n})^{d}}{(2r)^{d}}= \frac{1}{2}\frac{1}{\lambda_{n}^{d}} \xrightarrow[n\to \infty]{} \frac{1}{2}.
\end{split}
\end{equation*}

\noindent \emph{Case 2:} $r/\lambda_{n} > \alpha_{n}-d_{n}(x)$.

Now, the region $R_{n,i} \cap B_{r/\lambda_{n}}(x)$ is $A(r/\lambda_{n}, d_{n}(x), \alpha_{n})$ (see Figure \ref{FIG:example 3}). A direct computation of the derivative of $\Ha^{d}(A(r, a,c))$ with respect to $a$ (over the admissible range of values for $a$) shows that this quantity is minimized when $a = \max\{0, c-r\}$. Hence,
\begin{equation*}
\Ha^{d}(A(r/\lambda_{n}, d_{n}(x), \alpha_{n})) \geq \Ha^{d}(A(r/\lambda_{n}, \max\{0, \alpha_{n} - r/\lambda_{n}\}, \alpha_{n})).
\end{equation*}
If $r/\lambda_{n} \leq \alpha_{n}$, this region contains half of a $d$-dimensional ball of radius $r/\lambda_{n}$, so the same estimate as in Case 1 follows. Instead, if $r/\lambda_{n} > \alpha_{n}$ this region is $A(r/\lambda_{n}, 0, \alpha_{n})$, which consists of half of a $d$-dimensional ball of radius $r/\lambda_{n}$ minus a small spherical cap. In that case, using Eq. \ref{EQ:generalization inclusion 1} and part \textbf{(B)}, 
\begin{equation*}
\begin{split}
    \frac{\Ha^{d}\left(\Pi_{d} \cap B_{r}^{d+1}(\Pi_{d}(x))\right)}{(2r)^{d}} &\geq \frac{\Ha^{d}(A(r /\lambda_{n}, 0, \alpha_{n}))}{(2r)^{d}}\\
    & = \frac{\frac{1}{2}(2r/\lambda_{n})^{d}-\Ha^{d}(SC(r/\lambda_{n}, \alpha_{n}))}{(2r)^{d}}\\
    & \geq \frac{1}{2}\frac{1}{\lambda_{n}^{d}} - \frac{\Ha^{d}(SC(\alpha_{n}\lambda_{n-1}/\lambda_{n}, \alpha_{n}))}{(2r)^{d}} \xrightarrow[n \to \infty]{} \frac{1}{2},
\end{split}
\end{equation*}
because $\Ha^{d}(SC(\alpha_{n}\lambda_{n-1}/\lambda_{n}, \alpha_{n}))/(2r)^{d} \to 0$ when $n \to \infty$ as can be computed directly from Eq. \ref{EQ:area SC}.

\vsp
\noindent \textbf{(D)} $\ldd(\Pi_{d},\Pi_{d}(x)) \leq 1/2$ for $\Ha^{d}$-almost all $x \in [0,1]^{d}$.
\vsp

Let $x \in [0,1]^{d}$. For any $n \geq 1$, let $i$ be such that $x \in R_{n,i}$. If $y \in R_{n,i-1}\cup R_{n,i+1}$, using Eq. \ref{EQ:example lower bound} and $|x-y| \geq d_{n}(x)$, we obtain
\begin{equation*}
\begin{split}
   \dist(\Pi_{d}(x), \Pi_{d}(y)) & > \sqrt{d_{n}(x)^{2}+ \alpha_{n}^{2}\frac{1}{n^{2}}}  := r_{n}(x).
\end{split}
\end{equation*}
Taking this radius, we ensure
\begin{equation}
\label{EQ:generalization inclusion 2}
    \Pi_{d}\cap B^{d+1}_{r_{n}(x)}\left(\Pi_{d}(x)\right) \subset \Pi_{d}\left(R_{n,i} \cap B_{r_{n}(x)}(x)\right)
\end{equation}
Notice that the region $R_{n,i} \cap B_{r_{n}(x)}(x)$ is contained in a $d$-dimensional ball minus the spherical cap $SC(r_{n}(x), d_{n}(x))$. Therefore, using Eq. \ref{EQ:generalization inclusion 2} and part \textbf{(B)},
\begin{equation*}
\begin{split}
    \frac{\Ha^{d}\left(\Pi_{d} \cap B^{d+1}_{r_{n}(x)}\left(\Pi_{d}(x)\right)\right)}{(2r_{n}(x))^{d}} & \leq \frac{\Ha^{d}\left(\Pi_{d}(R_{n,i}\cap B_{r_{n}(x)}(x))\right)}{(2r_{n}(x))^{d}} \\
    & = \frac{\Ha^{d}\left(R_{n,i}\cap B_{r_{n}(x)}(x)\right)}{(2r_{n}(x))^{d}} \\
    & \leq 1 - \frac{\Ha^{d}\left(SC(r_{n}(x), d_{n}(x))\right)}{(2r_{n}(x))^{d}}.
\end{split}
\end{equation*}
Now, let $k$ be any positive integer. If $d_{n}(x) < \alpha_{n}/(kn)$, then substituting in Eq. \ref{EQ:area SC} we obtain
\begin{equation*}
\begin{split}
    \frac{\Ha^{d}(SC(r_{n}(x), d_{n}(x)))}{(2r_{n}(x))^{d}} \geq \frac{\Ha^{d}(SC(1, 1/k))}{2^{d}} := \beta_{k},
\end{split}
\end{equation*}
which tends to $1/2$ as $k \to \infty$ because $\Ha^{d}(SC(1, 1/k))$ tends to $2^{d-1}$, the $\Ha^{d}$ measure of half a $d$-dimensional ball with unit radius. In that case,
\begin{equation*}
\begin{split}
    \frac{\Ha^{d}\left(\Pi_{d} \cap B_{r_{n}(x)}(\Pi_{d}(x))\right)}{(2r_{n}(x))^{d}} \leq 1 - \beta_{k}.
\end{split}
\end{equation*}
Hence, we find that if $d_{n}(x) < \alpha_{n}/(nk)$ infinitely often, then $\ldd(\Pi_{d}, \Pi_{d}(x)) \leq 1- \beta_{k}$. Therefore, for any $\epsilon > 0$ there is $k \geq 1$ such that
\begin{equation*}
\begin{split}
    &\left\{x \in [0,1]^{d}\ :\ \ldd(\Pi_{d}, \Pi_{d}(x)) > \frac{1}{2}+\epsilon\right\}  \\
    \subset & \bigcup_{N \geq 1} \left\{x \in [0,1]^{d}\ :\ d_{n}(x) \geq \frac{\alpha_{n}}{nk}\text{ for all $n \geq N$}\right\}.
\end{split}
\end{equation*}
We will show next that the set of points $\mathcal{X}_{N} := \{x \in [0,1]^{d}\ :\ d_{n}(x) \geq \alpha_{n}/(nk)\text{ for all $n \geq N$}\}$ has zero $\Ha^{d}$-measure for any $k \geq 1$, so we finally obtain
\begin{equation*}
    \Ha^{d}(\{x \in [0,1]^{d} \ :\ \ldd(\Pi_{d}, \Pi_{d}(x)) > 1/2\}) = 0.
\end{equation*}
We have (for $N$ sufficiently large)
\begin{equation*}
    \mathcal{X}_{N} \subset  \bigcup_{i = 1}^{2^{N^{2}}-1}\bigcap_{m \geq N} A_{m}^{N,i},
\end{equation*}
where $A_{m}^{N,i}$, defined for $m \geq N$, is a subset of $[0,1]^{d}$ that we construct recursively as follows. We start with $A_{N}^{N,i} = R_{N,i}$. We remove the outer parts of $R_{N, i}$, which consist of the product of $[0,1]^{d-1}$ with the two outer subintervals of $I_{N, i}$ of length $\alpha_{N}/(k(N+1))$. Hence, we are keeping the interior part which consists of the product of $[0,1]^{d-1}$ with an interval of length $(1-2/(k(N+1)))\alpha_{N}$. Then we take $A_{N+1}^{N,i}$ to be the union of all the rectangles $R_{N+1,j}$ which are contained in that inner part. We repeat the process with each one of these new rectangles, and iterate $m-N$ times until we have constructed $A_{m}^{N,i}$. Notice that at each step,
\begin{equation*}
\Ha^{d}\left(A_{m}^{N,i}\right) \leq  \left(1-\frac{2}{km}\right)\Ha^{d}\left(A_{m-1}^{N,i}\right).
\end{equation*}
Consequently, by induction
\begin{equation*}
\Ha^{d}\left(A_{m}^{N,i}\right) \leq \Ha^{d}(R_{N, i}) \prod_{n = N+1}^{m}\left(1 - \frac{2}{kn}\right).
\end{equation*}
Since
\begin{equation*}
    \prod_{n = N+1}^{m} \left(1 - \frac{2}{nk}\right) \leq e^{-\frac{2}{k}\sum_{n = N+1}^{m}\frac{1}{n}} \xrightarrow[m \to \infty]{} 0,
\end{equation*}
we find $\Ha^{d}(A_{m}^{N,i}) \to 0$ as $m \to \infty$. Therefore, their intersection has zero $\Ha^{d}$-measure. It follows immediately that $\Ha^{d}(\mathcal{X}_{N}) = 0$.

\vsp

Finally, notice that parts \textbf{(C)} and \textbf{(D)} together prove that $\ldd(\Pi_{d}, \Pi_{d}(x)) = 1/2$ for $\Ha^{d}$-almost all $x \in [0,1]^{d}$, whereas strictly we require it to hold for $\Ha^{d}$-almost all points of $\Pi_{d}$. This is immediately granted by part \textbf{(A)}. Moreover, part \textbf{(B)} has the trivial consequence that $\Pi_{d}$ has finite and positive $\Ha^{d}$-measure, hence is a $d$-set.

\bibliographystyle{amsalpha}
\bibliography{bibliography.bib}

\end{document}